# Forecasting with Pairwise Gaussian Markov Models

Marc Escudier
SAMOVAR Telecom Sudparis
Institut Polytechnique de Paris
Palaiseau, France
marc.escudier@telecom-sudparis.eu

Ikram Abdelkefi
SAMOVAR Telecom Sudparis
Institut Polytechnique de Paris
Palaiseau, France
ikram.abdelkefi@telecom-sudparis.eu

Clément Fernandes
SAMOVAR Telecom Sudparis
Institut Polytechnique de Paris
Palaiseau, France
clement.fernandes@telecom-sudparis.eu

Wojciech Pieczynski
SAMOVAR Telecom Sudparis
Institut Polytechnique de Paris
Palaiseau, France
wojciech.pieczynski@telecom-sudparis.eu

*Abstract*—**Pairwise Markov Models (PMMs) extend the well-known Hidden Markov Models (HMMs). Being significantly more general, PMMs enable several types of processing, like Bayesian filtering or smoothing, similar to those used in HMMs. In this paper, we deal with Bayesian forecasting. The aim is to show analytically in the simple stationary Gaussian case that the extent results obtained with HMM can be improved. We complete contributions with a theoretical error study and two real examples we deal with. Experiments show that PMMs-based forecasting can significantly improve HMMs-based ones.**

## I. INTRODUCTION

The problem dealt with by this article is to forecast a future realization of a hidden sequence from past observations of an observed one. This general problem is recurrent in many fields such as question answering, financial forecasting, traffic congestion forecasting, or still weather forecasting. One key aspect of the problem is the modelization of the probabilistic links between hidden and observed data. Indeed, the modelization chosen must be rich enough, to best reflect the reality that we try to model, but also simple enough, to ensure that the computations associated with the estimation processes are feasible and relatively fast. One model that has been particularly popular in this context is the hidden Markov model (HMM). HMMs models are well-known and widely used in speech recognition, [1], [2], in image segmentation, [3], [4], in weather forecasting, [5], [6]. Versions of the continuous-time model have also been developed with numerous applications including decision processes, [7], [8], [9]. Despite their simplicity, HMMs turn out to be very robust and are sufficiently efficient in many situations.

Let us consider two stochastic sequences $X_{1:N} = (X_1, \ldots, X_N)$, $Y_{1:N} = (Y_1, \ldots, Y_N)$ taking their values in $\mathbb{R}$. To simplify notations, we assume them to be centered and with a variance of 1. Let us consider two models for the distribution of the sequence $Z_{1:N} = (Z_1, \ldots, Z_N)$ of couples $Z_n = (X_n, Y_n)$. $Z_{1:N}$ will also sometimes be written with $Z_{1:N} = (X_{1:N}, Y_{1:N})$. The first one is the well-known Hidden Markov model (HMM), and the second one is the Pairwise Markov Model (PMM). We assume both models are considered stationary, Gaussian, and Markovian. This means that their dependence graphs have shapes presented in Fig. 1, with $a = Cov[X_n, X_{n+1}]$, $b = Cov[X_n, Y_n]$, $c = Cov[Y_n, Y_{n+1}]$, $d = Cov[X_n, Y_{n+1}]$, and $e = Cov[X_{n+1}, Y_n]$, which are independent from $n = 1, \ldots, N$. We consider the problem of how to forecast $X_{n+k}$ from the observed $Y_{1:n}$.

PMMs extend HMMs, so that they may improve upon their results, [10]. In particular, they can significantly improve the quality of Bayesian segmentation, in the case of discrete hidden data, [11], [12], [13], [14], [15], [16], [17], [18], [19]. In continuous hidden data PMMs have been successfully applied in Kalman filtering, [20], [21], [22], [23], [24], [25], [26], [27], [28]. Finally, the only application of PMMs to forecasting is proposed in the case of discrete finite data [29]. PMM models having made a certain contribution to the resolution of segmentation and smoothing problems, it appears important to focus on their possible contribution to prediction.

The contribution of this article is to study the notion of prediction with PMMs, in the Gaussian linear case, and to make a comparison never studied in the literature with the classic HMC prediction model. For this, after having defined the model and the mathematical forms necessary for the prediction, we present a simulation study showing the relative interest of PMMs over HMMs, and we deal with real data examples of forecasting.

## II. FORECASTING WITH HMMS AND PMMS

Let's consider a stationary PMM $Z_{1:N} = (X_{1:N}, Y_{1:N})$. From the stationarity, the distribution of PMM is defined by the distributions of Gaussian $(X_1, Y_1, X_2, Y_2)^T$, which is defined by the covariance matrix, denoted with $\Gamma_{CMM}$.

$$\Gamma_{CMM} = \begin{bmatrix} 1 & b & a & d \\ b & 1 & e & c \\ a & e & 1 & b \\ d & c & b & 1 \end{bmatrix}. \tag{1}$$

HMM is a particular PMM, in that the covariance matrix $\Gamma_{HMM}$ defining its distribution is of the form (also see (8)-(10) for equivalent classic definition):

$$\Gamma_{HMM} = \begin{bmatrix} 1 & b & a & ab \\ b & 1 & ab & ab^2 \\ a & ab & 1 & b \\ ab & ab^2 & b & 1 \end{bmatrix}. \tag{2}$$

We see that PMMs are defined with 5 parameters, while HMMs are defined with two parameters. We can take any values for parameters $a, b, c, d,$ and $e$ since the matrix $\Gamma_{CMM}$ is positive definite.

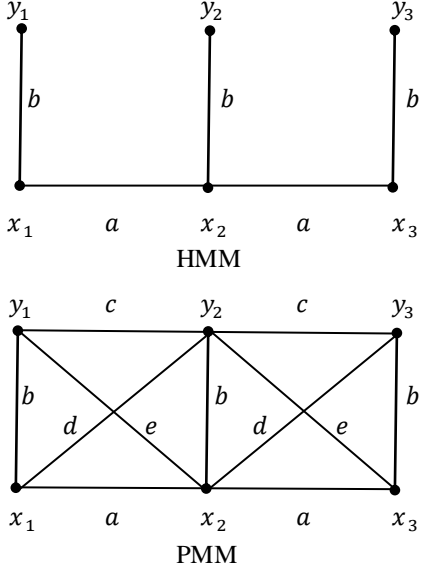

Fig. 1. Dependence graph of HMM and PMM. All means are 0, all variances are 1, $a, b, c, d,$ and $e$ designate covariances.

Classically, the PMM's distribution given with (1), is equivalently given in Markov form:

$$\begin{bmatrix} X_1 \\ Y_1 \end{bmatrix} \sim N\left(\begin{bmatrix} 0 \\ 0 \end{bmatrix}, \begin{bmatrix} 1 & b \\ b & 1 \end{bmatrix}\right); \tag{3}$$

$$\begin{bmatrix} X_{n+1} \\ Y_{n+1} \end{bmatrix} = A \begin{bmatrix} X_n \\ Y_n \end{bmatrix} + B \begin{bmatrix} U_{n+1} \\ V_{n+1} \end{bmatrix}, \tag{4}$$

with

$$A = \begin{bmatrix} a & e \\ d & c \end{bmatrix} \begin{bmatrix} 1 & b \\ b & 1 \end{bmatrix}^{-1}; \tag{5}$$

$$BB^T = \begin{bmatrix} 1 & b \\ b & 1 \end{bmatrix} - \begin{bmatrix} a & e \\ d & c \end{bmatrix} \begin{bmatrix} 1 & b \\ b & 1 \end{bmatrix}^{-1} \begin{bmatrix} a & d \\ e & c \end{bmatrix}; \tag{6}$$

and $U_1, V_1, \ldots, U_{N-1}, U_{N-1}$ normalized independent Gaussian variables.

As mentioned above, HMM is then obtained by setting

$$c = ab^2, \ d = ab, \ e = ab. \tag{7}$$

One can verify that by doing so we arrive at the classic representation of HMM:

$$X_1 \sim N(0, 1); \tag{8}$$

$$X_{n+1} = aX_n + \sqrt{1 - a^2} U_{n+1} \text{ for } n = 1, \ldots, N-1; \tag{9}$$

$$Y_{n+1} = bX_{n+1} + \sqrt{1 - b^2} V_{n+1} \text{ for } n = 0, \ldots, N-1. \tag{10}$$

Indeed, reporting (9) into (10) we have

$$Y_{n+1} = b(aX_n + \sqrt{1 - a^2} U_{n+1}) + \sqrt{1 - b^2} V_{n+1} =$$
$$abX_n + b\sqrt{1 - a^2} U_{n+1} + \sqrt{1 - b^2} V_{n+1}.$$

So that (9)-(10) can be written:

$$\begin{bmatrix} X_{n+1} \\ Y_{n+1} \end{bmatrix} = \tag{11}$$

$$\begin{bmatrix} a & 0 \\ ab & 0 \end{bmatrix} \begin{bmatrix} X_n \\ Y_n \end{bmatrix} + \begin{bmatrix} \sqrt{1-a^2} & 0 \\ b\sqrt{1-a^2} & \sqrt{1-b^2} \end{bmatrix} \begin{bmatrix} U_{n+1} \\ V_{n+1} \end{bmatrix},$$

which is of form (4) with

$$A = \begin{bmatrix} a & 0 \\ ab & 0 \end{bmatrix}; \tag{12}$$

$$B = \begin{bmatrix} \sqrt{1-a^2} & 0 \\ b\sqrt{1-a^2} & \sqrt{1-b^2} \end{bmatrix}; \tag{13}$$

$$BB^T = \begin{bmatrix} 1-a^2 & b(1-a^2) \\ b(1-a^2) & 1-a^2b^2 \end{bmatrix}. \tag{14}$$

We can verify that developing (6) while using (7) gives (14). To summarize, PMM is given with (3)-(4), while HMM is classically given with (8)-(10), which also is (3)-(4) with $c = ab^2$, $d = ab$, and $e = ab$.

Our aim is to compare HMMs with PMMs when applied to forecasting. More precisely, we forecast $X_{n+k}$ from $Y_{1:n}$. The forecasting is given with $E[X_{n+k}|Y_{1:n}]$, and its mean square error is given with $V[X_{n+k}|Y_{1:n}] = E[X_{n+k}^2|Y_{1:n}] - (E[X_{n+k}|Y_{1:n}])^2$. So, the problem is to compute $E[X_{n+k}|Y_{1:n}]$ and $E[X_{n+k}^2|Y_{1:n}]$, with HMM and PMM.

Let us specify the computation of $E[X_{n+k}|Y_{1:n}]$ in the case of PMM. We will use the formula

$$E[X_{n+k}|Y_{1:n}] = E[E[X_{n+k}|X_n, Y_{1:n}]|Y_{\{1:n\}}]. \tag{15}$$

To use (15), one computes $E[X_{n+k}|X_n, Y_{1:n}] = E[X_{n+k}|X_n, Y_n] = \varphi(X_n, Y_n)$, and then $E[\varphi(X_n, Y_n)|Y_{1:n}]$. Setting:

$$A^k = \begin{bmatrix} \alpha_{1,k} & \alpha_{2,k} \\ \alpha_{3,k} & \alpha_{4,k} \end{bmatrix}, \tag{16}$$

we have according to (4)

$$E[X_{n+k}|X_n, Y_n] = \alpha_{1,k} X_n + \alpha_{2,k} Y_n, \tag{17}$$

which comes from the fact that $U_{n+1}, V_{n+1}, \ldots, U_{n+k}, V_{n+k}$ are independent and centered Gaussian variables.

As $\varphi(X_n, Y_n)$ is linear in $X_n$, we have

$$E[\varphi(X_n, Y_n)|Y_{1:n}] = \varphi(E[X_n|Y_{1:n}]) = \tag{18}$$
$$\alpha_{1,k} E[X_n|Y_{1:n}] + \alpha_{2,k} Y_n.$$

As for $V[X_{n+k}|Y_{1:n}]$, it can be computed recursively with:

$$V[X_{n+1}, Y_{n+1}|Y_{1:n}] \text{ given}, V[X_{n+k+1}, Y_{n+k+1}|Y_{1:n}] = \tag{19}$$
$$AV[X_{n+k}, Y_{n+k}|Y_{1:n}]A^T + BB^T$$

## III. EXPERIMENTS

In this section, we first present a study comparing the theoretical error of a PMM and a HMM, with data following a PMM distribution. The goal of this study is to evaluate the potential gain, in terms of mean square error (MSE), from using a PMM model for forecasting, compared to an HMM model. Then we present two examples of real data forecasting, where using a PMM model is of interest compared to the HMM model.

### A. Theoretical error study

To illustrate the potential gain, in terms of MSE, from using a PMM model for forecasting, compared to an HMM model, we study the theoretical MSE of both models knowing that the data follow a PMM distribution. This particular choice seems to be of interest because it allows the computation of the theoretical error for both models; as the HMM model is a particular PMM, assuming that the data follow an HMM distribution would result in the same theoretical error. It is important to note that this study is biased toward PMM, as we suppose that the data follows a PMM distribution, so forecasting by PMM is the optimal solution. Furthermore, real data are not likely to follow one of the models. Nonetheless, it can give an idea of the interest of PMM in forecasting, as well as testing the robustness of HMM with respect to PMM, in the most favorable case for the PMM.

The theoretical MSE for PMC is simply $V[X_{n+k}|Y_{1:n}]$, which we can compute with the classical Kalman filtering for $k = 0$ and with (19) for $k \geq 1$. For HMM, one can compute the theoretical error knowing that the data follows a PMM distribution. Indeed, according to the projection theorem $E_{PMM}[X_{n+k}|Y_{1:n}]$ minimizes $E[\|X_{n+k} - V\|^2]$, for $V$ in $span(\{Y_1, ..., Y_n\})$. This is an orthogonal projection, which is characterized by the orthogonality of $X_{n+k} - E_{PMM}[X_{n+k}|Y_{1:n}]$ with each $U \in span(\{Y_1, ..., Y_n\})$. Noting that $E_{HMM}[X_{n+k}|Y_{1:n}] \in span(\{Y_1, ..., Y_n\})$, we can write according to Pythagoras' theorem:

$E[\|X_{n+k} - E_{HMM}[X_{n+k}|Y_{1:n}]\|^2] =$
$E[\|X_{n+k} - E_{PMM}[X_{n+k}|Y_{1:n}]\|^2] +$
$\quad E[\|E_{PMM}[X_{n+k}|Y_{1:n}] - E_{HMM}[X_{n+k}|Y_{1:n}]\|^2]$. (20)

Then $E[\|E_{PMM}[X_{n+k}|Y_{1:n}] - E_{HMM}[X_{n+k}|Y_{1:n}]\|^2]$ is computable in the case of stationary HMM and PMM with $X_1 \sim N(0,1), Y_1 \sim N(0,1)$, with a recursion on $n$, for $k = 0$, and an explicit formula in the function of $A^k$, for $k \geq 1$.

Let us set:

$BB^T = \begin{bmatrix} \beta_1 & \beta_2 \\ \beta_2 & \beta_3 \end{bmatrix}, A = \begin{bmatrix} \alpha_1 & \alpha_2 \\ \alpha_3 & \alpha_4 \end{bmatrix}$,

For $k = 0$, $E[\|E_{PMM}[X_n|Y_{1:n}] - E_{HMM}[X_n|Y_{1:n}]\|^2] =$

$\quad E[\|\sum_{i=1}^{n}(\alpha_{i,n}^{PMM} - \alpha_{i,n}^{HMM}) Y_i\|^2]$, (21)

with for all $1 \leq i \leq n-2$:

$\alpha_{i,n}^{PMM} = (\alpha_1 - \alpha_3 \frac{\alpha_1 \alpha_3 V[X_{n-1}|Y_{1:n-1}] + \beta_2}{\alpha_3^2 V[X_{n-1}|Y_{1:n-1}] + \beta_3}) \alpha_{i,n-1}^{PMM}$; (22)

$\alpha_{n-1,n}^{PMM} = \left(\alpha_1 - \alpha_3 \frac{\alpha_1 \alpha_3 V[X_{n-1}|Y_{1:n-1}] + \beta_2}{\alpha_3^2 V[X_{n-1}|Y_{1:n-1}] + \beta_3}\right) \alpha_{n-1,n-1}^{PMM} + \alpha_2 -$

$\quad \alpha_4 \frac{\alpha_1 \alpha_3 V[X_{n-1}|Y_{1:n-1}] + \beta_2}{\alpha_3^2 V[X_{n-1}|Y_{1:n-1}] + \beta_3}$; (23)

$\alpha_{n,n}^{PMM} = \frac{\alpha_1 \alpha_3 V[X_{n-1}|Y_{1:n-1}] + \beta_2}{\alpha_3^2 V[X_{n-1}|Y_{1:n-1}] + \beta_3}$. (24)

For $k \geq 1$, $E[\|E_{PMM}[X_{n+k}|Y_{1:n}] - E_{HMM}[X_{n+k}|Y_{1:n}]\|^2] =$

$\quad E[\|\sum_{i=1}^{n}(\alpha_{i,n+k}^{PMM} - \alpha_{i,n+k}^{HMM}) Y_i\|^2]$, (25)

with for all $1 \leq i \leq n-1$:

$\alpha_{i,n+k}^{PMM} = \alpha_{1,k} \alpha_{i,n}^{PMM}$; (26)

$\alpha_{n,n+k}^{PMM} = \alpha_{1,k} \alpha_{n,n}^{PMM} + \alpha_{2,k}$. (27)

The formulas (22)-(24) and (26)-(27) are also true respectively for $\alpha_{i,n}^{HMM}$ and $\alpha_{i,n+k}^{HMM}$, for $1 \leq i \leq n$, setting (12) and (14).

Finally, the covariances $E[Y_i Y_j]$, which intervene in (21) and (25), can be computed with:

For $1 \leq i \leq n$, $E[Y_i^2] = 1$;

$\quad$ For $i \neq j$, $E[Y_i Y_j] = b\alpha_{3,|i-j|} + \alpha_{4,|i-j|}$. (28)

There are two components, which can lead to a gain in terms of MSE for forecasting with PMM instead of HMM. The first component is the difference in computation of $E[X_{n+k}|X_n, Y_n]$. The second is the difference in computation of $E[X_n|Y_{1:n}]$, which acts as an initialization value for the prediction. To illustrate both effects, we present the variation of both theoretical MSE in the function of $n$, with $k = 0$, in Fig. 2 and Fig. 4. The variation of both theoretical MSE in the function of $k$, with different values of $n$, are presented in Fig. 3 and Fig.5. The parameters are chosen in the following manner. For the first batch of experiments (Fig. 2 and Fig. 3), we construct an HMM with $a = 0.90$, $b = -0.20$, so $c, d$ and $e$ are fixed to $c = ab^2$, $d = e = ab$. Then we construct a PMC by changing the value of e with $e = ab - 0.4$, while $d$, and $c$ stay the same. For the second batch of experiments (Fig. 4 and Fig. 5), we begin with the same HMC with $a = 0.90$, $b = -0.20$, and then construct a new PMC with $e = ab - 0.4$ and $d = ab - 0.2$, c still valued at $ab^2$.

The first observation we can make is that filtering with a PMM can improve significantly upon HMM filtering. Indeed, in Fig. 2 the PMM MSE can become up to six times smaller than the one of HMM as $n$ increases, and in Fig. 4 the PMM MSE is eventually more than ten times smaller than the one of HMM. If we look at the evolution of the MSE in function of $k$, when $n = 1$ (blue and black curves of Fig. 3 and Fig. 5) we can see that the PMM can improve upon forecasting with HMM, but the effect is one to two magnitudes smaller than for the filtering, with a maximum gain of around ten percent on the MSE.

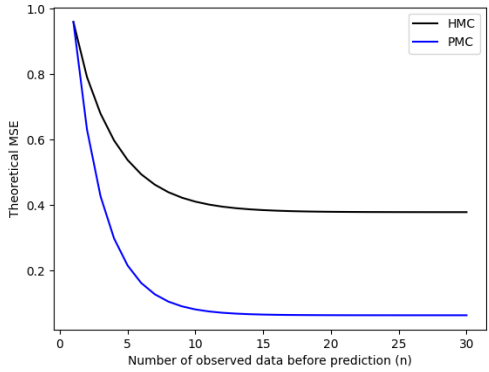

Fig. 2. Theoretical mean square error of filtering for the HMM and PMM models, with parameters $a = 0.90, b = -0.20, c = ab^2, d = ab, e = ab - 0.4$, for data following the PMC distribution, in the function of the number of observations $(n)$.

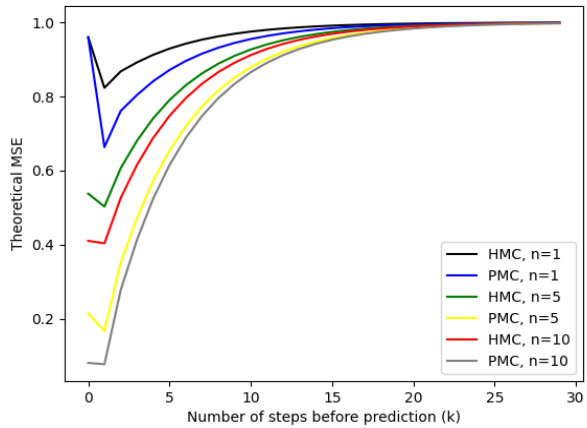

Fig. 3. Theoretical mean square error of prediction for the HMM and PMM models, with parameters $a = 0.90, b = -0.20, c = ab^2, d = ab, e = ab - 0.4$, for data following the PMM distribution, in function of the number of steps before prediction $(k)$, with different numbers of observation $n = 1, 5, 10$.

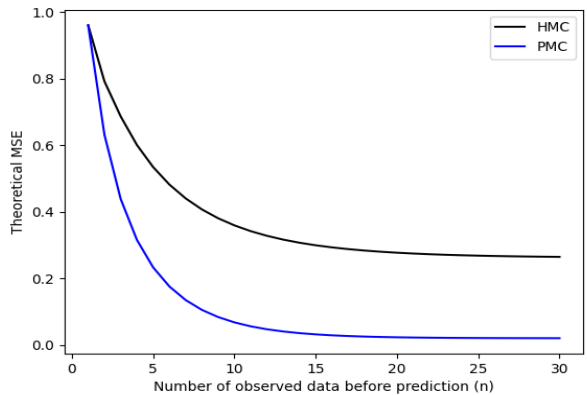

Fig. 4. Theoretical mean square error of filtering for the HMM and PMM models, with parameters $a = 0.90, b = -0.20, c = ab^2, d = ab - 0.2, e = ab - 0.4$, for data following the PMM distribution, in function of the number of observations $(n)$.

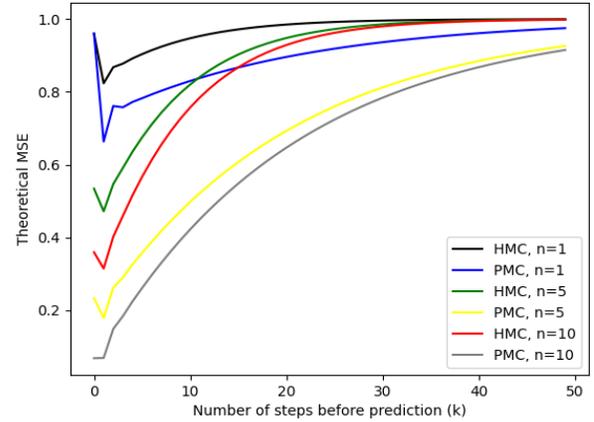

Fig. 5. Theoretical mean square error of prediction for the HMM and PMM models, with parameters $a = 0.90, b = -0.20, c = ab^2, d = ab - 0.2, e = ab - 0.4$, for data following the PMM distribution, in function of the number of steps before prediction $(k)$, with different numbers of observation $n = 1, 5, 10$.

However, if we look at Fig. 5, we can see that the convergence to the MSE of the common stationary marginals of both models when $k$ increases is slower for the PMM model, which might be interesting when forecasting at distant horizons. Finally, we can observe the effect of both factors on the MSE with the green, red-yellow, and gray curves of Fig. 3 and Fig. 5. In these cases, the filtering component greatly improves upon the predictions of the PMM compared to the HMM ones for the first few values of $k$, but this gain decreases as $k$ increases as expected. The most interesting point is in Fig. 5, where the gain from the filtering component combined with the slower convergence of the MSE of PMM to the MSE of the marginal law leads to an MSE for the PMM that remains significantly smaller than the MSE of the HMM, even for large values of $k$. To conclude, in the most favorable case for the PMM, the latter can improve both the filtering and prediction with HMM. However, the gain in the filtering phase is much more significant than the one in the prediction phase. Moreover, in some cases, the PMM MSE converges slower to the MSE of the stationary law than the HMM one, when $k$ increases. This last effect, combined with the gain in the filtering phase can make the prediction of PMM significantly better than the one of HMM, at distant horizons, for some specific correlation values.

### B. Forecasting real data

We propose some forecasting experiments using historical weather data for simple univariate variables. The aim is to give two examples of prediction where the classical HMM is outperformed by the PMM described above.

HMM are commonly used models to compute time series forecasting, [30], with applications in many fields: epidemiology, [31], load forecasting, [32], and weather forecasting, focusing on different variables like temperature, wind speed, or photovoltaic energy, [33], [34], [35]. Extreme weather, climate, and water-related events caused 11,778 reported disasters between 1970 and 2021, with just over 2 million deaths and US$ 4.3 trillion in economic losses, according to the WMO, [36]. In this context of climate change, two important variables can be considered: atmospheric pressure and soil moisture. Accurate forecasting of these two variables is of crucial importance to climatology, agriculture, or meteorology. Atmospheric pressure is a good indicator of extreme weather events like tornadoes, hurricanes, and extreme precipitation. Soil moisture is a good indicator of drought which has consequences in irrigation planning, water demand, and, more widely, on social and economic environment. Soil moisture seems also to be a source of seasonal predictability in temperate regions, [37]. Therefore, many studies proposing effective models including Decision Trees, Random Forest, Neural networks, and more particularly LSTMs, have been carried out in order to predict these indicators, [38], [39], [40], [41]. More precisely, the prediction problems including the atmospheric pressure and especially the soil humidity are often studied with HMM models, [42], and, in particular, with Kalman models such as the Ensemble Kalman Filter, [43], the Unscented Kalman Filter, [44], or the Extended Kalman Filter, [45]. All these models derive from the classical Kalman model based on HMC. It therefore appears interesting and useful to test and compare the prediction of these two variables with the extension that we propose, the PMM model, in order to observe the behavior of our model in this task. This also seems useful because, by showing the advantage of the PMM over the HMM, even more complex Kalman models cited above could be improved by substituting the assumptions of the HMM with those of the PMM which are less restrictive.

We present a simple application of the stationary PMM and HMM models for the prediction of the two variables mentioned before by temperature, based on an open-source meteorological dataset issued by the site Meteo Blue. It consists of hourly values of these three parameters for a station based in Basel (Switzerland) taken over two years, from January 2021 to January 2023 (17544 hourly values). The application is based on two experiences: the first one is the prediction of the atmospheric pressure by the temperature and the second one is the prediction of the soil moisture by the temperature. So, it means that the observed variable $Y_n$ is the temperature in the two experiments, and the hidden variable we want to predict at the horizon $k$, $X_{n+k}$, is the atmospheric pressure in the first one and the soil moisture in the second one.

To quantify the non-stationarity of our time series, we first applied a Dickey-Fuller test: we obtained a p-value of 0 for both atmospheric pressure and soil humidity series and a p-value of 0.05 for the variable $Y$ describing temperature (we can reject the hypothesis of non-stationarity if the $-value \leq 0.05$). To prevent periodic temperature variations from polluting the forecast, we chose to remove the cyclical components (day-night oscillations, seasons) from the variable $Y$. The seasonal component is estimated using the result of the particular linear regression $f$ of $Y$ given by $f(i) = \theta_0 + \theta_1 \cos\frac{i2\pi}{24} + \theta_2 \sin\frac{i2\pi}{24} + \theta_3 \cos\frac{i2\pi}{8772} + \theta_4 \sin\frac{i2\pi}{8772}$, for $i$ in $[1,17544]$ obtained with the least squares method. We find $\theta = (\theta_0, \theta_1, \theta_2, \theta_3, \theta_4)$ minimizing the sum of the squares of the residuals by writing $\theta = (J^T W J)^{-1} J^T W Y$ where:

$$J = \begin{pmatrix} \cos(\frac{1 \times 2\pi}{24}) & \cdots & \sin(\frac{1 \times 2\pi}{8772}) & 1 \\ \vdots & \ddots & \vdots & \vdots \\ \cos(\frac{N \times 2\pi}{24}) & \cdots & \sin(\frac{N \times 2\pi}{8772}) & 1 \end{pmatrix},$$

with dimensions of $(N \times 5)$ with $N = 17544$, and

$$W = \begin{pmatrix} 1/\sigma & \cdots & 0 \\ \vdots & \ddots & \vdots \\ 0 & \cdots & 1/\sigma \end{pmatrix},$$

with dimensions of $(N \times N)$ and $\sigma$ is the empirical variance of $Y$. By considering $Y - f$ instead of $Y$, the daily oscillations are contained in the prediction phase. Furthermore, a second Dickey-Fuller test was carried out after removing the seasonal component: the p-value obtained may indicate that with this transformation (we obtained a p-value of 0 for $Y - f$), the model has come closer to our theoretical stationary case.

Then, to stay close to our theoretical study, the data are standardized. We estimate the parameters $a, b, c, d, e$ on a chosen month of the first year of our dataset, starting from the 1st January 2021 at 00:00 (672 hourly values for February and 744 or 720 for the others months) with the classical empirical estimators of variance and covariance. Next, we use them to define our two models, HMM and PMM, as described before and compute the predictions $\hat{X}_{n+k} = E[X_{n+k}|Y_{1:n}]$ of values (function of the values $Y_1, Y_2, \ldots, Y_n$ and of $a, b, c, d, e$) on the same month of the second year. Finally, to evaluate the quality of the prediction, we calculate the associated mean square error (MSE)

$$\frac{1}{M}\sum_{i=1}^{M}|\hat{X}_i - X_i|^2$$

on the standardized values (standardized MSE), with $M$ equal to $24 \times L$ and $L$ representing the number of days of the chosen month. The results of the experiments will be presented for the month of January. So, we have $L = 31, M = 744$; $a, b, c, d, e$ are estimated on the month of January 2021 and the prediction is tested on the month of January 2022.

For the first experiment, we obtained the following values for the parameters: $a = 0.996, b = -0.6, c = 0.986, d = -0.599, e = -0.602$. We can notice that the link between $X_i$ and $Y_i$ doesn't exactly correspond to a linear correlation for $X_{i+1}$ and $Y_i$, and as for $X_i$ and $Y_{i+1}$. Besides, we can verify that the data clearly does not correspond to the HMM condition since we have: $a \times b^2 = 0.35$ and $c = 0.986$. We present the MSE values in TABLE I. for different values of chain size $n$ and different values of prediction horizon $k$. It seems that PMM outperforms HMM in all these cases. In addition, the advantage of forecasting with PMM rather than with HMM appears more significant with a large number of observations n and a more distant horizon k. Increasing the horizon k and varying n seem to have less impact on the performance of the PMM compared to that of the HMM. We also present two graphs of predicted values for the fourth week of January (second year): the first Fig. 6, corresponds to a prediction with a lag of two days (48 hours), and the second in Fig. 7, corresponds to a prediction with a lag of one day (24 hours).

On these graphs, we see that the forecasts with the PMM are slightly better than with the HMM. The HMM model appears more sensitive to data variations by oscillating more strongly than the PMM.

For the second experiment, the prediction of soil moisture, we obtained the following values for the parameters: $a = 0.996, b = 0.543, c = 0.986, d = 0.545, e = 0.542$. In this case, too, we can see that we are far from the HMM model conditions since $a \times b^2 = 0.29$ and $c = 0.986$. The MSE values for this experiment are presented in TABLE II. As in the first experiment, we can make the observation that forecasting with the PMM is better than with the HMM. In this case, the variations of n and k have a minor effect on the performance of the two models. Surprisingly, predictions at a distant horizon seem better than those at a nearer horizon in some cases. As in the first experiment, we present two graphs of predicted values for the second week of January: the first Fig. 8, corresponds to a prediction with a lag of three days (72 hours), and the second Fig. 9, with a lag of two days (48 hours). As in the first experiment, the HMM model oscillates more strongly than the PMM, moreover, it has a tendency to underestimate the values, while the PMM estimations are closer to the real values.

TABLE I. STANDARDIZED MSE FOR THE PREDICTION OF ATMOSPHERIC PRESSURE (IN HECTOPASCAL) IN FUNCTION OF TEMPERATURE (IN CELSIUS) FOR THE MONTH OF JANUARY-SECOND YEAR

| n | k | MSE HMC | MSE PMC |
|---|---|---------|---------|
| 5 | 10 | 0.85 | 0.59 |
| 5 | 24 | 0.77 | 0.57 |
| 5 | 48 | 1.15 | 0.83 |
| 20 | 10 | 0.75 | 0.53 |
| 20 | 24 | 0.82 | 0.57 |
| 20 | 48 | 1.24 | 0.78 |
| 50 | 10 | 0.69 | 0.59 |
| 50 | 24 | 0.86 | 0.59 |
| 50 | 48 | 1.15 | 0.64 |

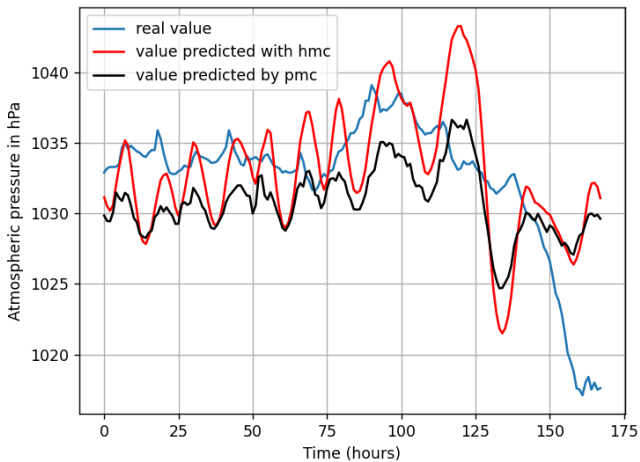

Fig. 6. Prediction of atmospheric pressure by temperature at horizon $k = 48$ with $n = 5$ for the 4[th] week. Standardized MSE for HMC: 0.25; standardized MSE for PMC: 0.21.

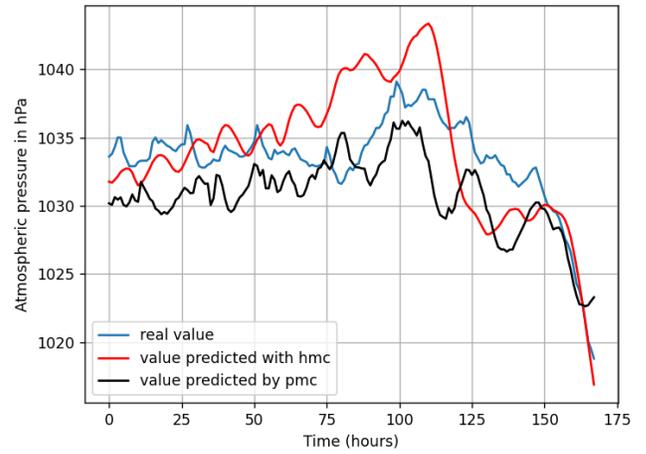

Fig. 7. Prediction of atmospheric pressure by temperature at horizon $k = 24$ with $n = 20$ for the 4[th] week. Standardized MSE for HMC: 0.12; standardized MSE for PMC: 0.11.

TABLE II. STANDARDIZED MSE FOR THE PREDICTION OF SOIL MOISTURE (RATIO OF THE VOLUME OF WATER CONTAINED IN A SUBSTRATE AND THE VOLUME OF THE SUBSTRATE) IN FUNCTION OF TEMPERATURE (IN CELSIUS) FOR THE MONTH OF JANUARY-SECOND YEAR

| n | k | MSE HMC | MSE PMC |
|---|---|---------|---------|
| 5 | 10 | 1.31 | 0.79 |
| 5 | 24 | 1.02 | 0.70 |
| 5 | 48 | 1.05 | 0.77 |
| 20 | 24 | 1.13 | 0.67 |
| 20 | 48 | 1.22 | 0.74 |
| 20 | 72 | 1.20 | 0.69 |
| 50 | 24 | 1.05 | 0.68 |
| 50 | 48 | 1.17 | 0.64 |
| 50 | 72 | 1.05 | 0.60 |

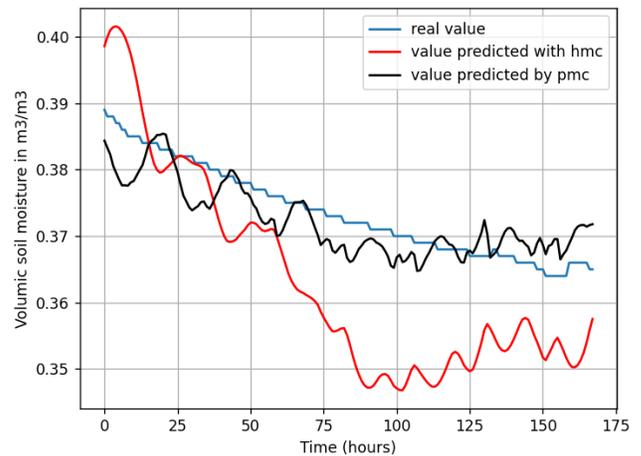

Figure 8. Prediction of soil moisture by temperature at horizon $k = 72$ with $n = 20$ for the 2[nd] week. Standardized MSE for HMC: 0.44; standardized MSE for PMC: 0.03.

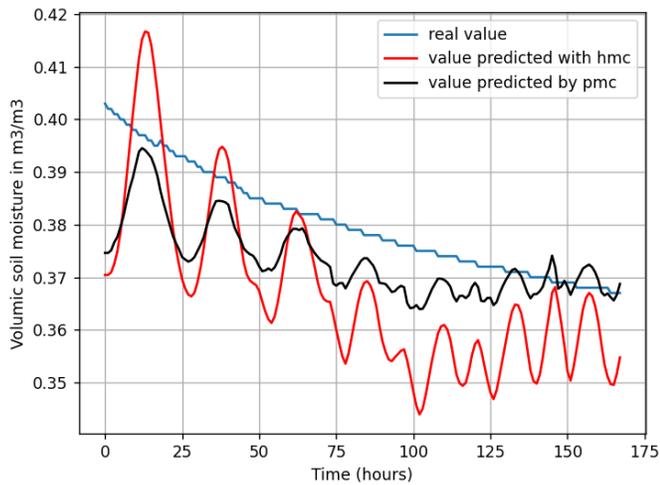

Fig. 9. Prediction of soil moisture by temperature at horizon $k = 48$ with $n = 5$ for the 2$^{nd}$ week. Standardized MSE for HMC: 0.65; standardized MSE for PMC: 0.21.

IV. CONCLUSIONS AND PERSPECTIVES

In this paper, we dealt with forecasting with pairwise Gaussian Markov models and with hidden ones. We first studied the theoretical error of both models for data following a PMM distribution. While these studies were biased toward the PMM model, they nonetheless showed that PMM-based forecasting methods can significantly improve upon HMM-based ones, while also being equivalent to the latter in the worst case. We then compared the two models on two real applications, the forecasting of atmospheric pressure and soil moisture, which are crucial variables in climatology, agriculture, or meteorology, by temperature. The results confirmed our theoretical conclusion as in both cases the PMM-based forecasting is significantly better than the HMM-based one. It is important to note that characteristics of our model remain to be studied; in particular the exact influence of variations in the size of the chain and the prediction horizon on the quality of the prediction. It can also be noted that the model presented is deliberately simple and is therefore not suitable for the study of numerous real data. However, the presented comparisons will probably be useful to understand the contribution of PMMs in more complex Kalman models.

As a perspective, we may consider more complex switching models, well adapted to nonstationary data, [46], [47], and non-Gaussian data based on Kalman models, [48]. Indeed, to our knowledge, such models are not commonly used for forecasting, but they are likely to further improve the methods described in this article.